\documentclass[a4paper,10pt,twoside]{article}

\usepackage[utf8]{inputenc}

\usepackage[english]{babel}
\usepackage{indentfirst}
\usepackage{graphicx}
\usepackage{amssymb,amsmath}
\usepackage{multicol}

\usepackage{a4wide}
\usepackage{subfig}
\usepackage{amsfonts}
\usepackage{tabulary}
\usepackage{latexsym}
\usepackage{amsthm}

\usepackage{amsfonts}
\usepackage{mathrsfs}
\usepackage{url}
\usepackage{caption}

\begin{document}
	
\title{Banking risk as an epidemiological model:\\ an optimal control approach\thanks{This 
work is part of first author's Ph.D., which is carried out at University of Aveiro
under the Doctoral Program in Applied Mathematics MAP-PDMA, 
of Universities of Minho, Aveiro and Porto.}
\thanks{This is a preprint of a paper whose final and definite form is in
\emph{Operational Research}, Springer Proceedings in Mathematics \& Statistics, 
available at [http://www.springer.com/series/10533]. 
Paper Submitted 23/March/17; Revised 29/May/17; Accepted 11/July/2017.}}
	
\author{Olena Kostylenko$^{1}$\\
{\tt \small o.kostylenko@gmail.com}
\and Helena Sofia Rodrigues$^{1,2}$\\
{\tt \small sofiarodrigues@esce.ipvc.pt}
\and Delfim F. M. Torres$^1$\\
{\tt \small delfim@ua.pt}}


\date{$^1${\text{Center for Research and Development in Mathematics and Applications (CIDMA)}},
Department of Mathematics, University of Aveiro, 3810-193 Aveiro, Portugal\\[0.1cm]
$^2$School of Business Studies, Polytechnic Institute of Viana do Castelo,\\
4930-678 Valença, Portugal}

\maketitle


\begin{abstract}
\noindent The process of contagiousness spread modelling is well-known in epidemiology. 
However, the application of spread modelling to banking market is quite recent. In this work, 
we present a system of ordinary differential equations, simulating data from the largest European banks. 
Then, an optimal control problem is formulated in order to study the impact of a possible 
measure of the Central Bank in the economy. The proposed approach enables qualitative 
specifications of contagion in banking obtainment and an adequate analysis and prognosis 
within the financial sector development and macroeconomic as a whole. We show that 
our model describes well the reality of the largest European banks. Simulations 
were done using \textsf{MATLAB} and \textsf{BOCOP} optimal control solver, 
and the main results are taken for three distinct scenarios.\\
			
\end{abstract}		

		
\section{Introduction}

Mathematical models of spread of epidemics are widely used in different fields of studies 
and directly related to everyone's life. They are well used in medicine for the purpose 
of studying an epidemiological process, analysis and forecasting outbreaks of such infections 
as influenza, cholera, HIV/AIDS, syphilis and others \cite{[1],[2],[3],{MR3476198}}; 
in computer science, for studying the spread of computer viruses \cite{[5],[4]}; 
in psychology, for studying how individuals could change 
their behaviour in response to others (for example, an applause, 
contagious yawning, the use of social networks such as MySpace, Facebook and others) 
\cite{[6],[7]}; in marketing, for modelling a viral distribution of advertising 
\cite{[8],[9],[10]}; in economics and finance, in order to study the financial 
crisis transmission, contagion in banking \cite{[11],[13],[12]}; and even 
in science fiction, to describe the propagation of zombies attacks \cite{[15],[14],[16]}.
		
The economic field has gained a special spotlight in recent years. This is due 
to the global economic decline during the early 21st century. In terms of overall impact, 
it was the worst global downturn since Second World War, by the conclusions of the 
International Monetary Fund. The crisis touched every country and it has led to the financial 
sector's collapse in the world economy, the effects of which can be seen 
and felt till today \cite{[19],[17],[18]}.
		
Several scientists, from all over the world, have proposed their theories 
about the development of the financial crisis and methods to prevent it. 
However, this topic has not been studied yet completely. Indeed, there 
is a lack of consensus among scientists, and financial crises continue 
to occur from time to time. These arguments support the importance, relevance 
and necessity to continue to study such phenomena.
		
This paper is focused on studies of contagion in banking market using 
an epidemiological approach through mathematical modelling. Initially, 
it is important to understand that the contemporary economy is an open system, 
which is based on direct and inverse, vertical and horizontal linkages, 
and can be developed successfully only with effective management 
of these relations at both macro and micro levels \cite{[20]}.
		
Banks serve the needs of market participants, and constitute the horizontal 
linkages in the market environment. They do not just form their own resources, 
but also provide domestic accumulation of funds for the national economy development. 
Therefore, banks organize the monetary-lending process (moving money from lenders 
to borrowers) and issue banknotes. Banks are thus, constantly, faced with many 
risks in their activity. As an example, the credit risk (risk of not fulfilling 
liabilities to the credit institution by a third party) may occur during the 
loan or other equivalent actions, which are reflected on the balance sheet 
and could have the off-balance character. The occurrence of this risk can cause 
infection and even bank failure. That could trigger the contagion of other banks 
and the emergence of a banking crisis.
		
The risk of contagion in the financial sector presents a serious threat 
for a country's economy. Failure of one particular institution in the 
financial sector threatens the stability of many other institutions. 
This situation is called a systemic risk \cite{[23]}. Moreover, 
the spread of one type of crisis is able to initiate the development 
of other types of crises as well. For example, the external debt crisis 
could undermine the stability of banks; problems in the banking sector 
can launch a debt crisis; and the banking crisis, in its turn, often 
precedes currency crisis. Thus, the risky lending and loan defaults, 
generally, precede crisis in banking \cite{[22]}. Banking crises often 
are accompanied by panic among population. It usually happens when 
many banks suffer runs simultaneously, as people have suspicions 
and mistrust to banks, and they suddenly try to withdraw 
their money from bank accounts.
		
Banks are directly at the centre of the financial system and it is easy 
to understand that a crisis in banking is one of the most serious 
type of financial crisis. Therefore, it is very important to study 
how the crisis spreads in banking market, to find its basic laws 
and methods to control it. A good understanding of their propagation 
mechanism makes possible to find and propose suitable policy interventions, 
which can most effectively reduce their contagious spread. This is the 
main goal of our research. Namely, we show that the application 
of epidemiological models is able to describe well the nature 
and character of the contagion spread and its behaviour over time. 
Our analysis identifies which measures to adopt and when they must 
be taken in order to prevent effects and serious negative consequences 
for a particular bank, and for economy as a whole.
		
The paper is structured as follows. In Section~\ref{s2}, 
we introduce the main bank concepts into the epidemic model. 
Various scenarios of contagion and the results of simulations  
are presented in Section~\ref{s3}. An optimal control problem
is then formulated in Section~\ref{s4}, where the control 
simulates the European Central Bank. In Section~\ref{s5}, 
the main conclusions are carried out.


\section{The SIR mathematical model}
\label{s2}

The availability of interactions between banks in a financial system
confirms the possibility of a contagion occurrence. 
Contagion, by definition, refers to the idea that any type of financial crisis
may spread from one subject (financial institution) to another. 
As an example, if a large number of bank customers are withdrawing their funds,
then it provokes a bank run that may spread from a few banks to many others 
within the country, or even spreading from one country to another. 
Currency crisis is one type of financial crisis, which also may spread 
from one country to another. Even sovereign defaults and stock market crashes 
spread across countries. Such processes of contagion in economy are very similar 
to the disease propagation in a population (from one individual to another). 
This similarity allows us to consider contagion in such type of financial sectors, 
as banking, using the same mathematical models of infection spreading as used in epidemiology.

There are many kinds of epidemiological models of spread of infection, 
which differ in the different assumptions about the disease's nature and character, 
about the structure of population that are under consideration, and others constraints 
in which the model is based. In our work, we assume that contagion can be transmitted
from an infected bank to another one, which has not yet been infected, and after recovery 
the bank produces immunity for a long time. Our assumptions about the process of contagion 
spreading in the banking sector are very similar to the characteristics of many childhood diseases, 
including chickenpox, measles, rubella and others, after which a strong immunity is produced. 
In order to simulate an epidemiological process of such diseases, there is a well-known suitable model, 
called SIR, which was first introduced in the works of Kermack and McKendrick \cite{[21]}. 
The model is named SIR for the reason that population is divided into three classes:
class $S$ of susceptible, those who are susceptible to the disease; class $I$ of infected,
those who are ill/infected; class $R$ of recovered, those who are recovered 
and have become immune to the disease.

Mathematically, the SIR model is a multi-parameter system of ordinary differential equations,
where each equation defines the current state of health of an investigated object.
This model is relevant and frequently used in current research.
Our work is also based on its application to the transmission of a contagious disease between banks.
The model must be, however, adapted to the banking sector, 
in order to have a deeper understanding of how the process of infection transmission 
is carried out, and especially how banks, as financial institutions, become contagious.

The earlier formulated rules of bank transitions from one condition to another
lead to the following system of differential equations:
\begin{equation}
\label{eqSIR}
\begin{cases}
\displaystyle \frac{dS(t)}{dt} = - \beta S(t) I(t),\\[0.3cm]
\displaystyle \frac{dI(t)}{dt} = \beta S(t) I(t) - \gamma I(t),\\[0.3cm]
\displaystyle \frac{dR(t)}{dt} = \gamma I(t),
\end{cases}
\end{equation}
$t \in [0,T]$, subject to the initial conditions
\begin{equation}
\label{eqIC}
S(0) = S_0, \quad I(0) = I_0, \quad R(0) = R_0.
\end{equation}
The first equation defines the number of banks that leave the group $S$ 
by becoming a target of market speculation at time $t$; the second equation 
defines the number of contagious banks at time $t$; and the third equation 
defines the number of banks that have recovered from the crisis at time $t$, 
where $t\in [0,T]$, $T > 0$.

Model \eqref{eqSIR}--\eqref{eqIC} depends on two parameters. 
Parameter $\beta$ denotes the contagion spreading rate. 
It represents the strength of contagion, where a susceptible bank, 
that is in contact with a single infected bank, changes its state 
to infected with probability $\beta$. The second parameter, $\gamma$, 
denotes the speed of recovery. It represents the resistance to contagion. 
Therefore, an infected bank changes its state to recover with probability $\gamma$.

A susceptible bank from group $S$ can obtain contagion if it has a relationship 
with a contagious bank from group $I$, and if it has not enough money in reserve 
to cover possible risk losses. In other words, the bank may be contaminated, 
if it has not enough strong ``health'' in order to resist an epidemic. 
According to the standard SIR methodology, we are assuming that all banks are similar, 
and both parameters $\beta$ and $\gamma$ are constant for the entire sample. However, 
in reality, these data are unique for each bank and depend on the force of contagion 
of the affected bank and the financial stability of the susceptible. Dynamics of the 
strength of banks is not taken into consideration in the present study. 

It is assumed that $N$, the number of institutions during the period of time under study,
is fixed throughout the contamination time (in our simulations, $N=169$ --- see Section~\ref{s3}). 
It means that $S(t)+I(t)+R(t) = N$ for all $t\in [0, T]$. Moreover, we assume that the initial 
conditions \eqref{eqIC} at the beginning of the period under study satisfy
\begin{equation}
\label{eq:IC:assump}
S(0)\gg I(0)>0 \quad \text{ and } \quad R(0)=0.
\end{equation}


\section{SIR model simulations with real data}
\label{s3}

The contagion and recovery propagation rates used in our paper, $\beta$ and $\gamma$, 
follow the statistics and empirical findings of Philippas, Koutelidakis and Leontitsis \cite{[24]}. 
The statistical data is taken, with respect to the year 2012, for the 169 largest European banks 
with total assets over than 10 bn Euro. The values of parameters $\beta$ and $\gamma$  
were calculated in \cite{[24]} by assuming that all banks of a country are infected; 
then, using Monte Carlo simulations, the parameters $\beta$ and $\gamma$ were tracked 
for each bank in a country; and aggregated in such a way they represent the average over simulations, 
with respect to the country in which the bank originates. Therefore, using the information 
that the total number of banking institutions ($N$) is equal to 169, and assuming that 
only one bank is contagious at the initial time, the values of initial conditions 
for the SIR model were obtained: $S(0)=168$, $I(0)=1$, and $R(0)=0$. The parameters 
$\beta$ and $\gamma$ are presented for each bank with respect to the country 
in which the bank originates.

In our work, in order to present contrast scenarios of contagion spreading among European banks,
we have chosen three countries that have completely 
different values of parameters $\beta$ and $\gamma$: 
\begin{enumerate}
\item in the first scenario, the initially contagious bank is located in Portugal;
	
\item in the second scenario, the group $I$ of infected  
consists only of a bank from Spain;
	
\item while in the third scenario the contamination begins 
from a United Kingdom's bank. 
\end{enumerate}
Such countries were chosen from \cite{[24]}. However, it should
be mentioned that our scenarios are different from those of \cite{[24]}. 
Indeed, in \cite{[24]} they begin with one or more random banks among three different levels of assets. 
In contrast, here we do not look to the assets level of the banks and, instead, we chose countries 
by their contagion and recovery rates, respectively $\beta$ and $\gamma$. 
Moreover, in our cases all the three scenarios begin with exactly 
one random infected bank in the chosen country, while in \cite{[24]} more than one is allowed.

\begin{figure}[ht]
\centering
\includegraphics[scale=0.8]{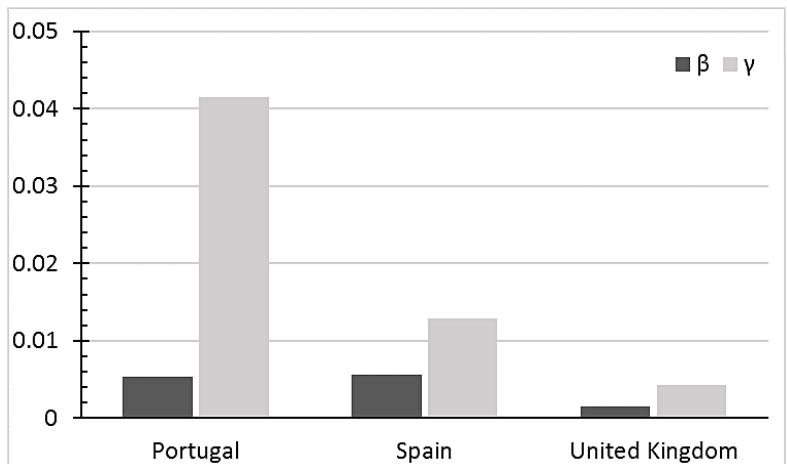}
\vspace*{0.3cm}
\caption{Summary statistics for $\beta$ and $\gamma$ based on the data of \cite{[24]}.}
\label{f1}
\end{figure}

\begin{figure}[ht]
\centering
\subfloat[Portugal]{\label{f2a}\includegraphics[width=0.33\textwidth]{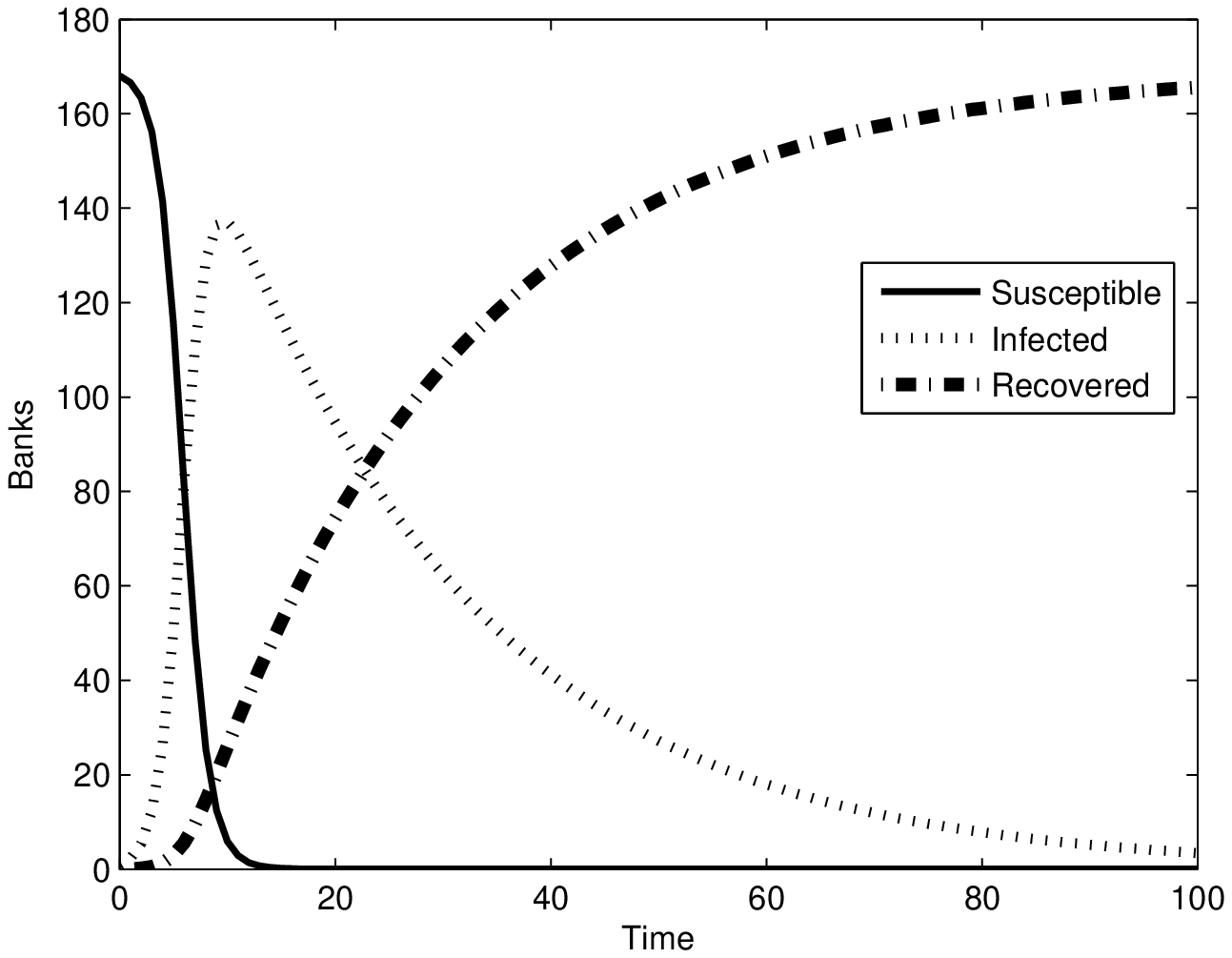}}
\subfloat[Spain]{\label{f2b}\includegraphics[width=0.33\textwidth]{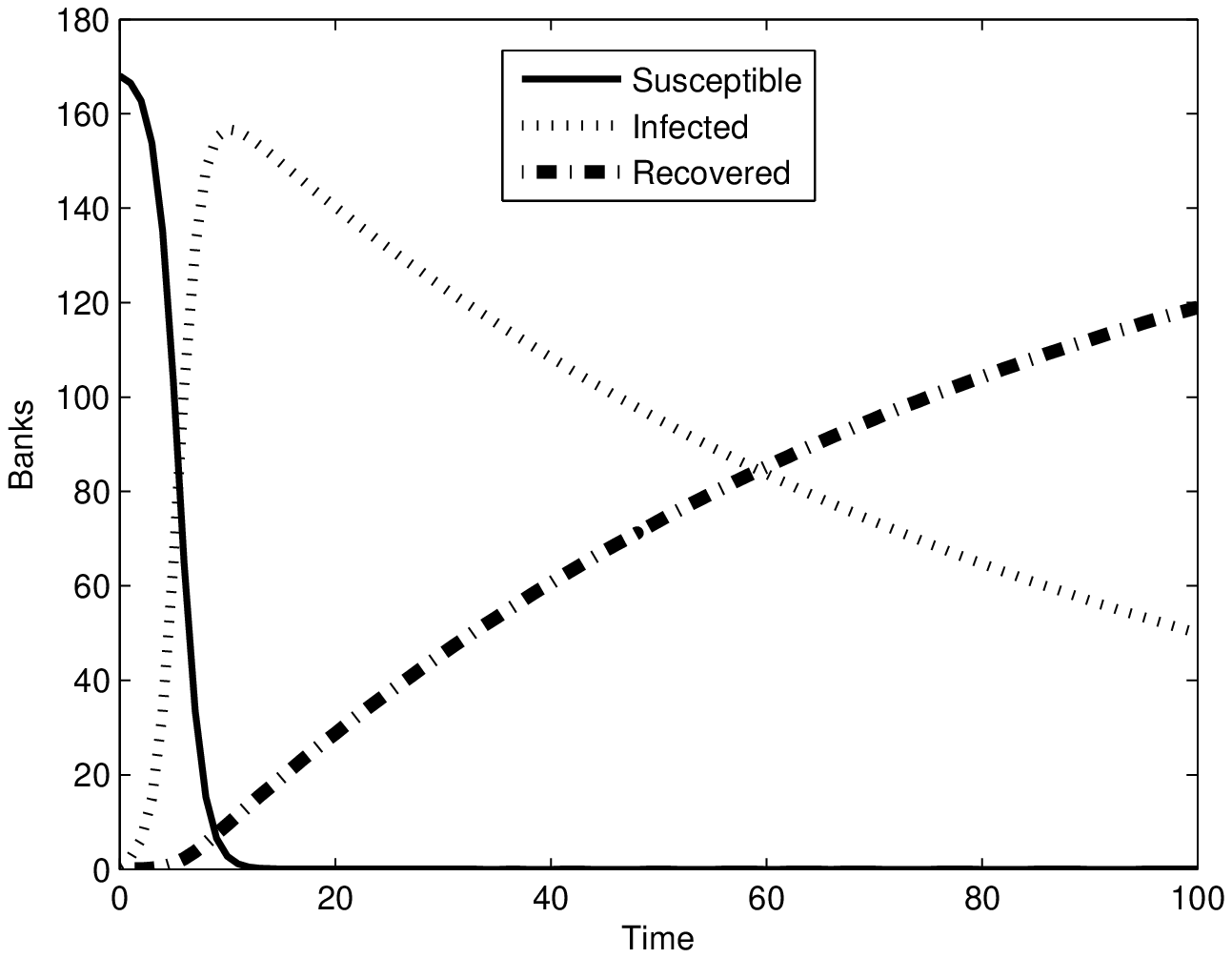}}
\subfloat[UK]{\label{f2c}\includegraphics[width=0.33\textwidth]{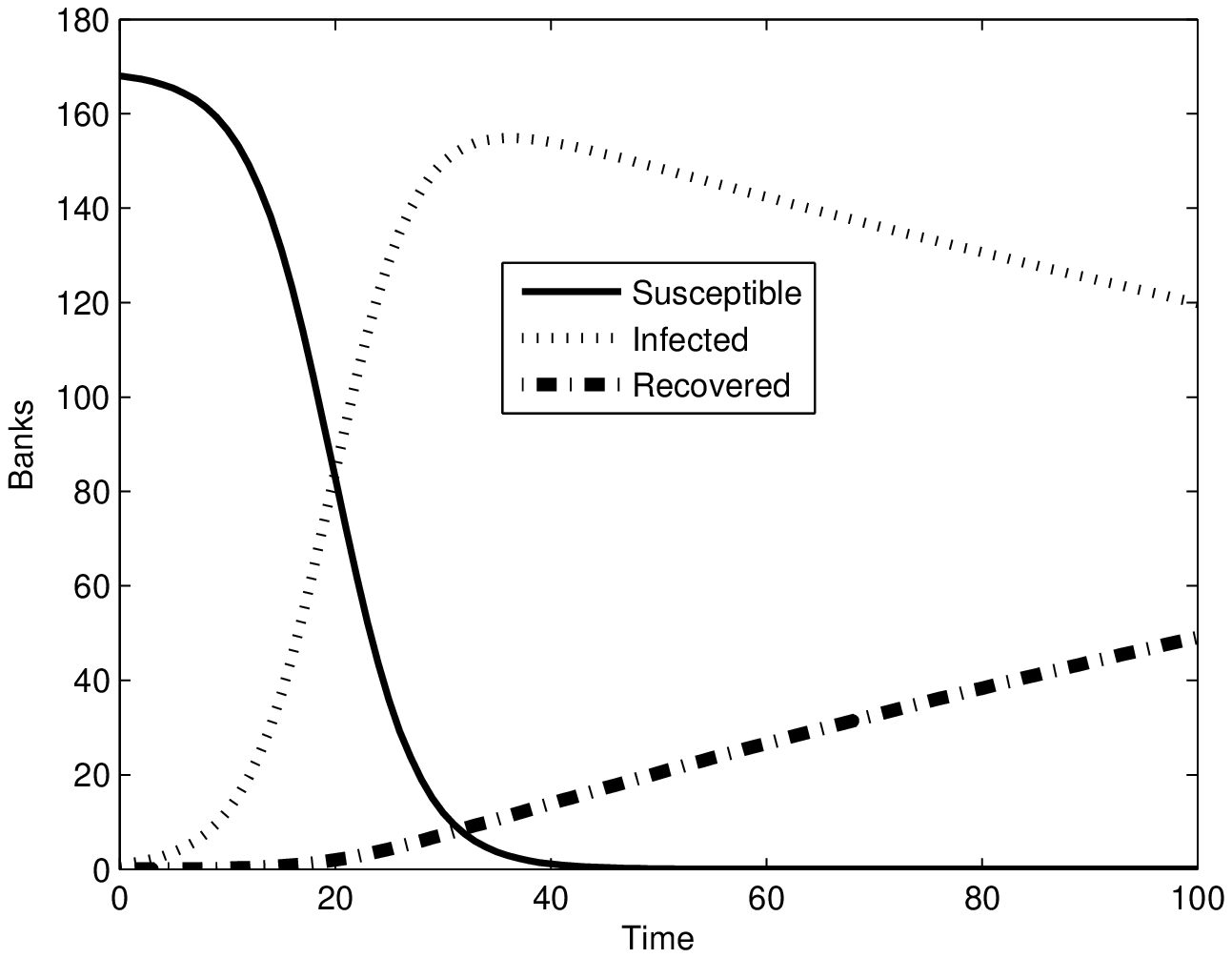}}
\vspace*{0.3cm}
\caption{The SIR contagion risk model \eqref{eqSIR}--\eqref{eqIC} with 
parameters $\beta$ and $\gamma$	as in Figure~\ref{f1}, $S(0)=168$, 
$I(0)=1$, $R(0)=0$ and $T=100$.}
\label{f2}
\end{figure}

Data of the contagion spreading rate and the speed of recovery 
for these three scenarios are shown in Figure~\ref{f1}
and the behaviour of the contagion risk model \eqref{eqSIR}--\eqref{eqIC},
during the same time scale of $T=100$ days, is shown in Figure~\ref{f2}.

In the first scenario, contagion spreads rapidly. 
In a short span of time, contagion reaches its peak and then swiftly decreases. 
The recovery process takes fast. See Figure~\ref{f2a}.
Regarding the second scenario, Figure~\ref{f2b} demonstrates that contamination occurs 
rapidly as in the first scenario. The contagion affects a large number of banks 
but the recovery process takes longer. In Figure~\ref{f2c}, we see what happens in the third
scenario: the spread of contagion occurs in a longer period of time and, consequently, 
the recovery process goes slowly too. The graphs of Figure~\ref{f2} show that 
the processes of contagion take place in different ways, depending on the country where it begins.

\begin{figure}[ht]
\centering
\subfloat[Portugal, $T=365$]{\label{f3a}\includegraphics[width=0.33\textwidth]{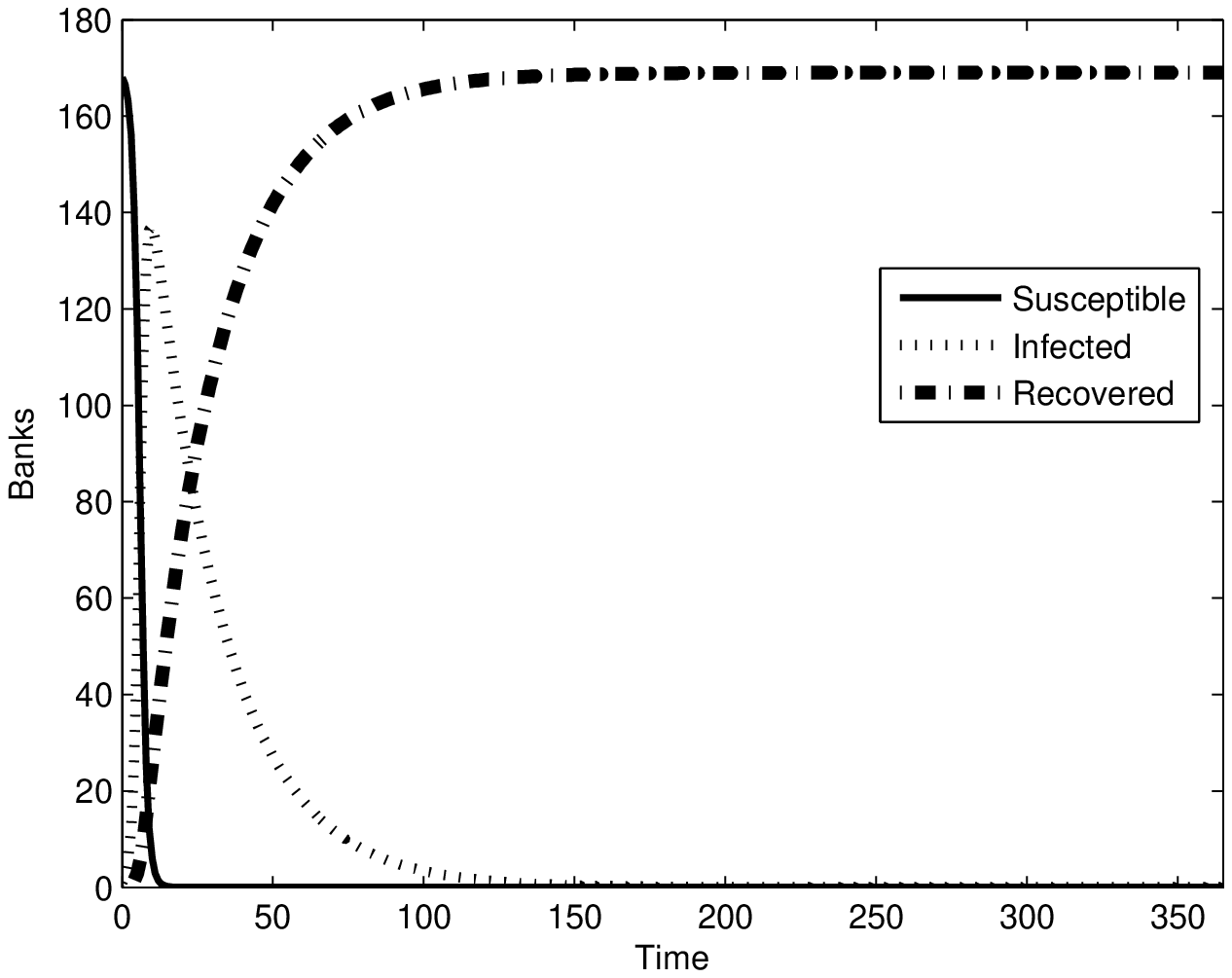}}
\subfloat[Spain, $T=450$]{\label{f3b}\includegraphics[width=0.33\textwidth]{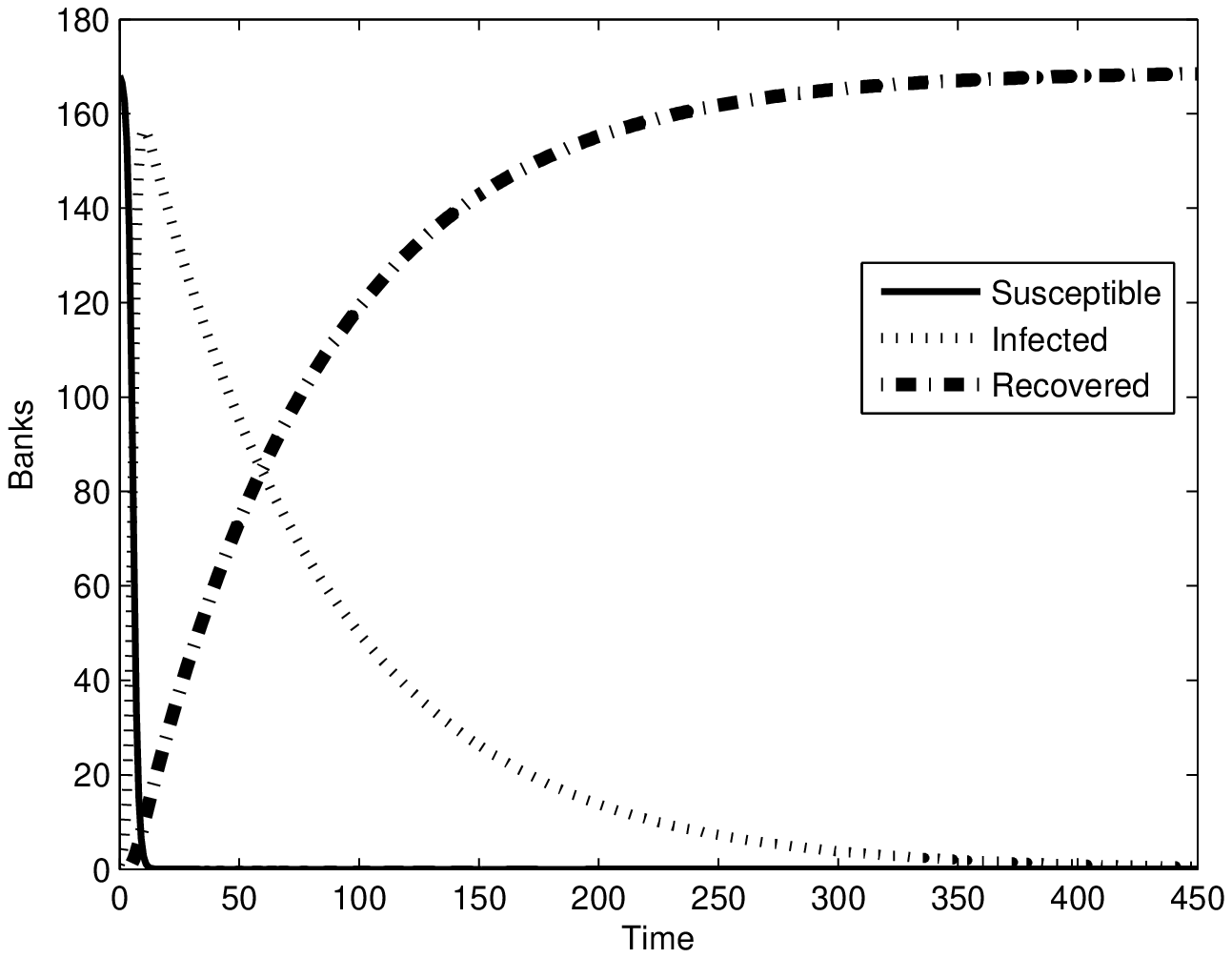}}
\subfloat[UK, $T=1200$]{\label{f3c}\includegraphics[width=0.33\textwidth]{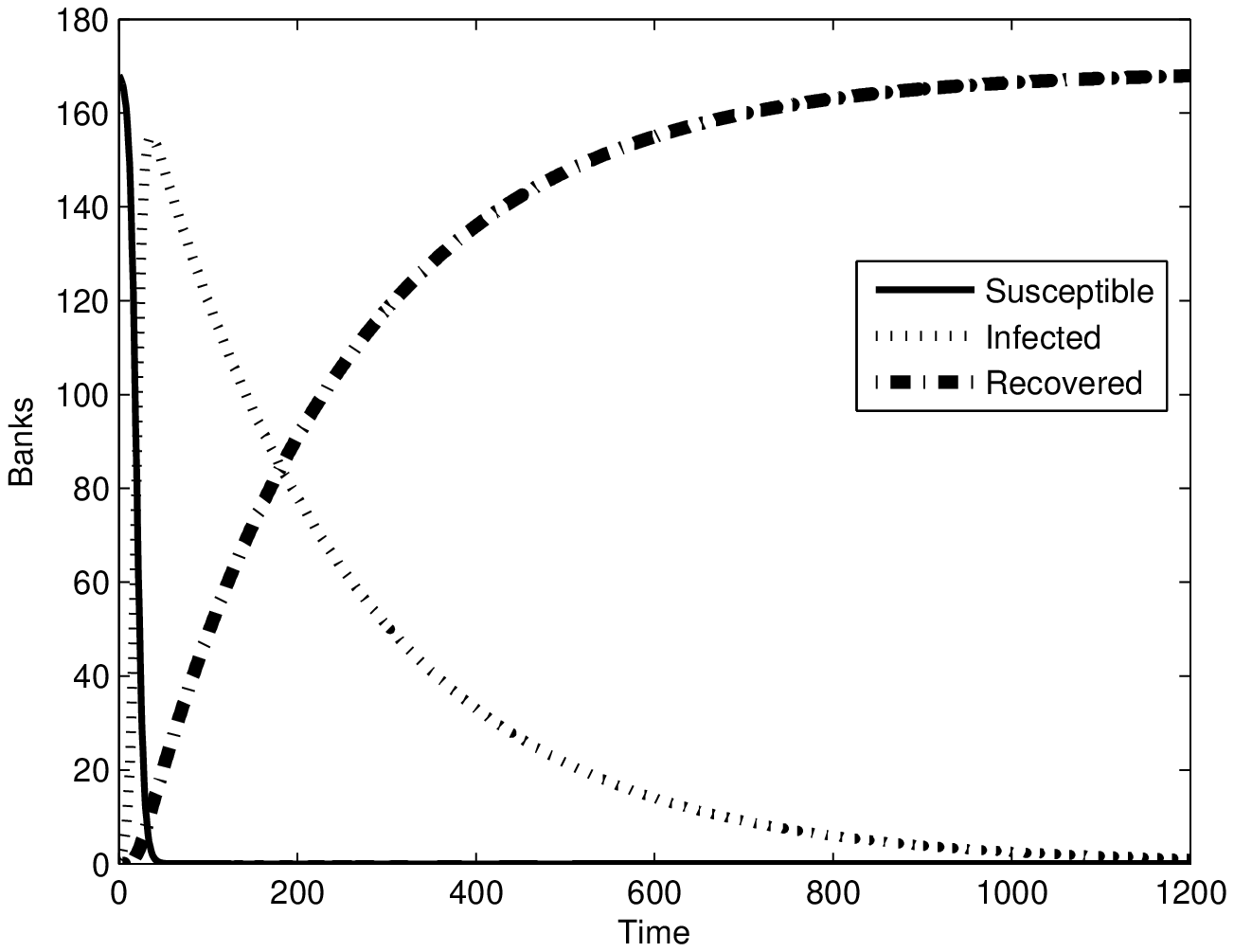}}
\vspace*{0.3cm}
\caption{The SIR contagion risk model \eqref{eqSIR}--\eqref{eqIC} with 
parameters $\beta$ and $\gamma$	as in Figure~\ref{f1}, $S(0)=168$, 
$I(0)=1$ and $R(0)=0$: illustration of the time needed to achieve 
the contagion-free equilibrium.}
\label{f3}
\end{figure}

Figure~\ref{f2} also shows that, in case of the first scenario, contagion 
has almost reached the contagion-free equilibrium ($I(T)=0$) at the end of 100 days. 
However, $T=100$ is not enough to reach the contagion-free equilibrium  
for the second and third scenarios. In those cases, it makes sense to increase the time. 
The results of such simulations are reflected in Figure~\ref{f3}.

As seen from Figure~\ref{f3a}, in case of bank contagion starting in Portugal,
less than in half a year the contagion spreading stops. Figure~\ref{f3b} 
demonstrates that a spreading of contagion starting in Spain will only stop after a full year. 
The third scenario is the most severe: the contagion will disappear and the banks will be recovered 
only after three years (see Figure~\ref{f3c}).

The reason behind the differences found, for the three scenarios considered,
rely in the different economic relevance, in the global banking market,
of the countries where contagion begins. So, if one of the large and important 
banks of United Kingdom will be the starting point in the spread of contagion, 
then the global banking market will experience serious problems, the negative 
effects of which will continue during a long period of time. In other words,
United Kingdom is one of the main financial markets to have into account
at European level. On the other hand, Portugal has a less powerful 
and influential economy on an European scale.


\section{Optimal control}
\label{s4}

Given the horizontal linkages in market environment, all banks are equally needed 
to be under financial supervision. Such supervision is necessary in order to prevent 
global contamination and avoid serious consequences due to spread of contagion 
between banks. This partially explains why the Central European Bank exists. 
The Central Bank is a supervisory competent authority that implements vertical 
connections among interacting commercial banks. Since one of the main goals 
is to avoid the wide dissemination of a contagion, we pose
the following optimal control problem:
\begin{equation}
\label{eqopt1}
\min \mathcal{J}[I(\cdot),u(\cdot)]
=I(T)+\int\limits_0^T b u^2(t) dt \longrightarrow \min
\end{equation}
subject to the control system
\begin{equation}
\label{eqopt2}
\begin{cases}
\displaystyle \frac{dS(t)}{dt} = - \beta S(t) I(t),\\[0.3cm]
\displaystyle \frac{dI(t)}{dt} = \beta S(t) I(t) - \gamma I(t)-u(t)I(t),\\[0.3cm]
\displaystyle \frac{dR(t)}{dt} = \gamma I(t)+u(t)I(t),
\end{cases}
\end{equation}
where $S(0) = S_0$, $I(0) = I_0$, $R(0) = R_0$,
satisfying \eqref{eq:IC:assump}, and the weight $b > 0$, are fixed.
In our simulations, the initial values of the state variables 
are given as discussed in Section~\ref{s3}, 
while the weight $b$, associated with the cost of control measures, 
is taken to be $1.5$, motivated by the value of possible recapitalization 
with state funds considered in \cite{[24]}. 
The first term of the objective functional \eqref{eqopt1} reflects 
the fact that the control organization should take care about the 
number of contagious institutions $I(T)$ at the final time $T$. 
The integral represents the general cost of financial assistance necessary 
to prevent the spread of contagion and economic decline in the period $[0, T]$. 
The control $u(\cdot)$ is a Lebesgue function with values in the compact set 
$[0, 1]$: $0 \leq u(t) \leq 1$, with $u(t)$ denoting the rate at which assistance 
will be provided to contagious banks, that is, it is the ratio between the financial 
support from the Central Bank at time $t$ and the financial needed by the banks at that time. 
In this way, $u(t) = 1$ means full support from the Central Bank at time $t$
(all money needed by the banks is being covered by Central Bank), while $u(t) = 0$ 
means no financial lending or recapitalization from the Central Bank at time $t$.

\begin{figure}[ht]
\centering
\subfloat[without control]{\label{f4a}\includegraphics[width=0.33\textwidth]{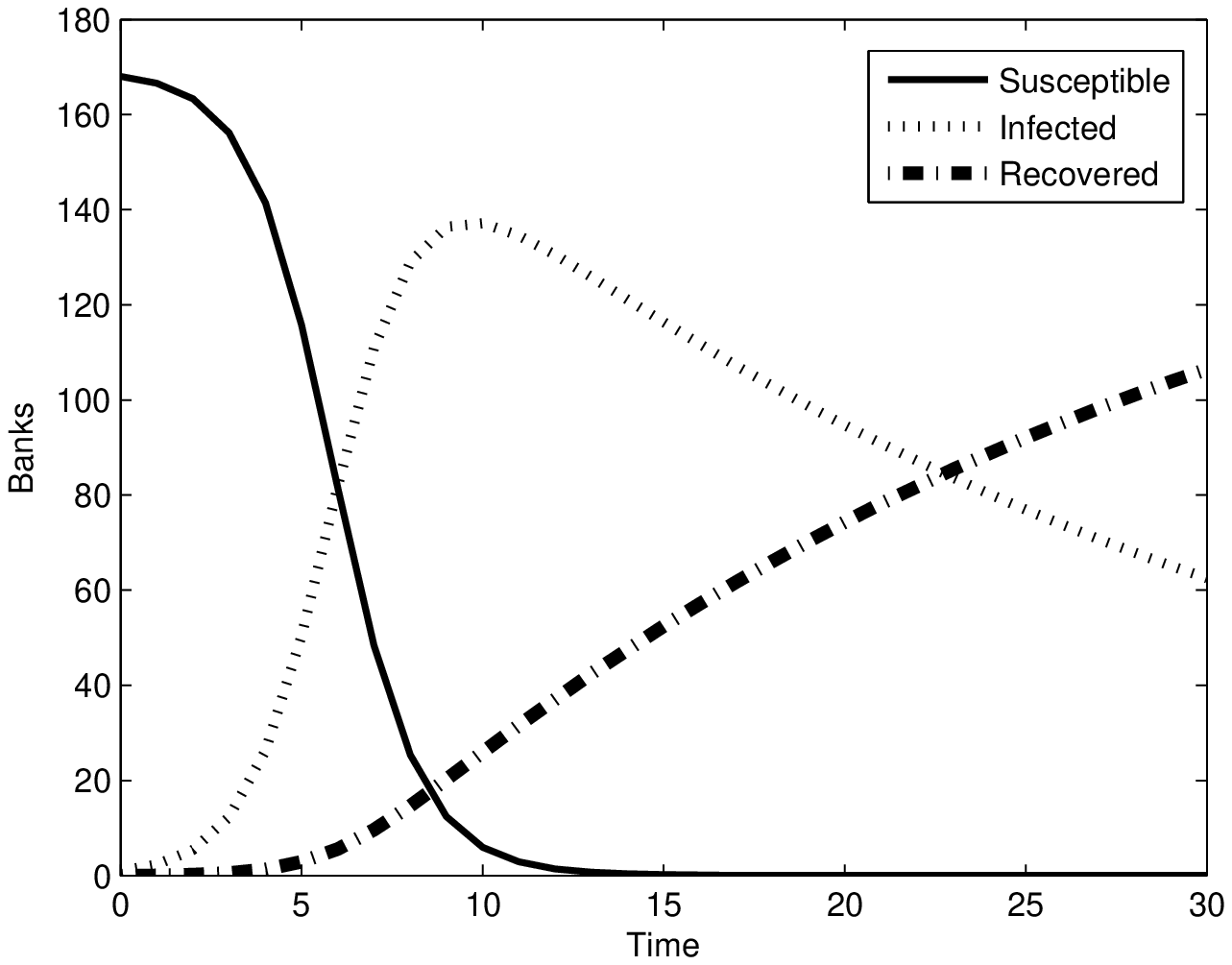}}
\subfloat[with optimal control]{\label{f4b}\includegraphics[width=0.33\textwidth]{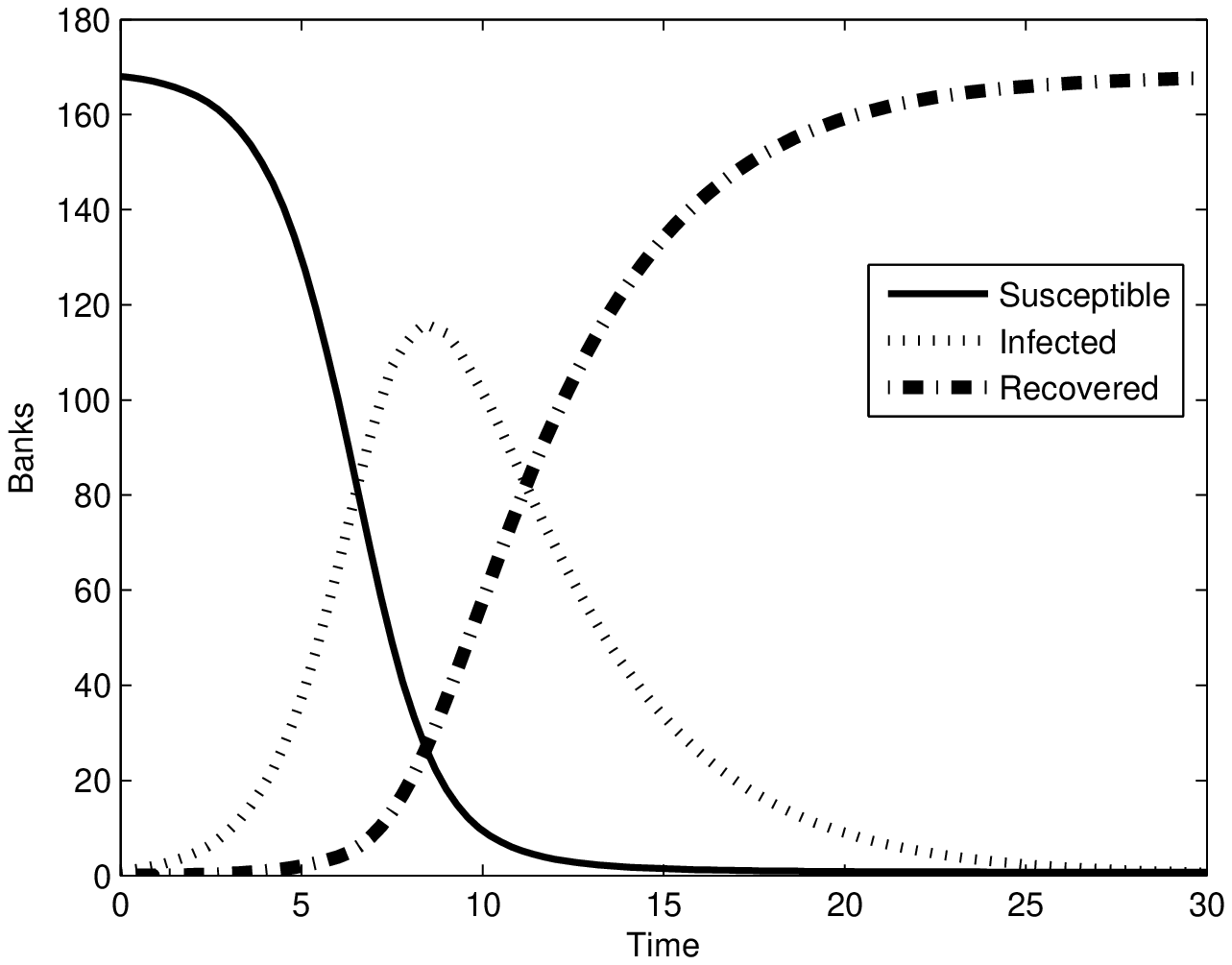}}
\subfloat[the extremal control $u(t)$]{\label{f4c}\includegraphics[width=0.33\textwidth]{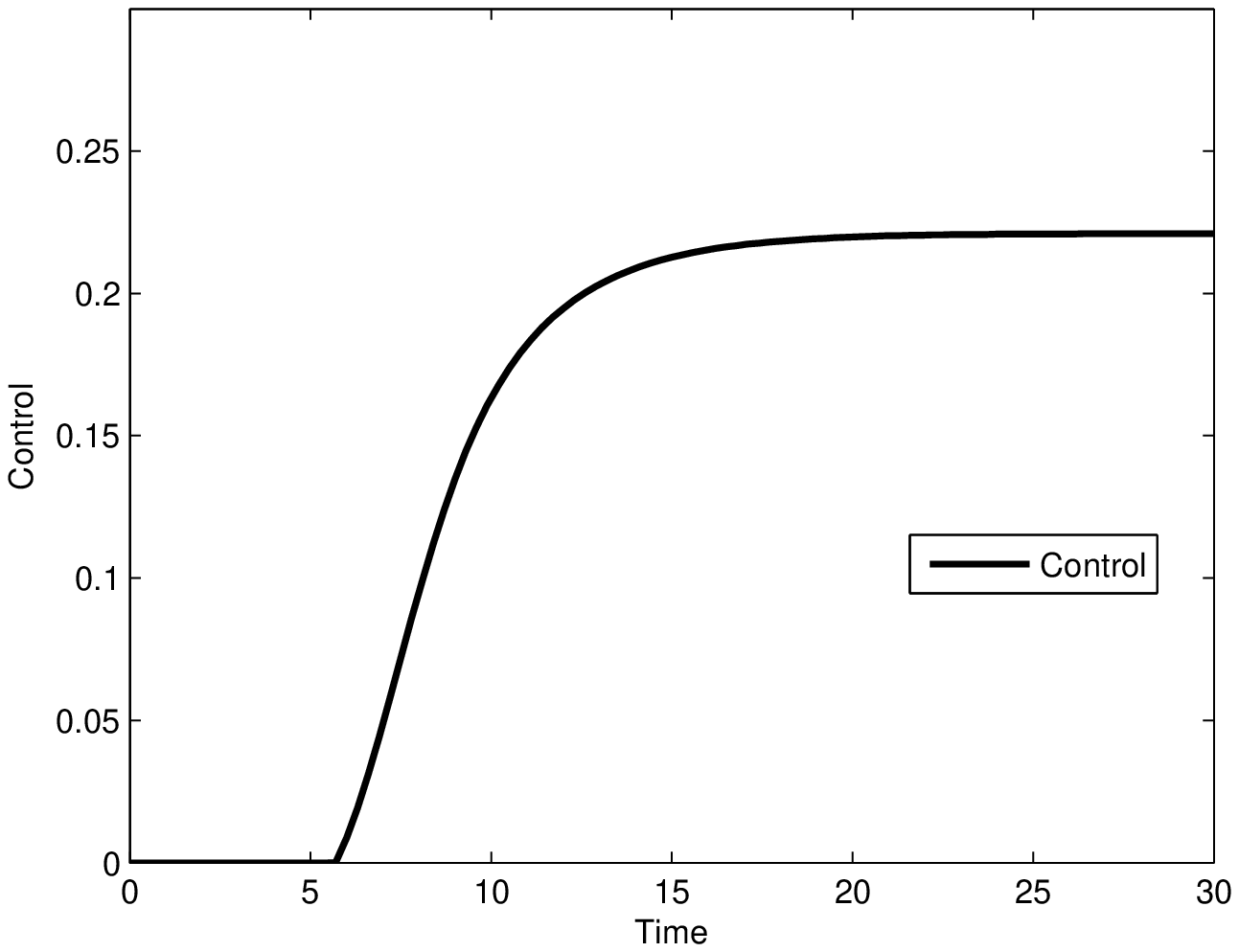}}
\vspace*{0.3cm}
\caption{Contagion risk from Portugal (Scenario~1) with and without optimal control.}
\label{f4}
\end{figure}

\begin{figure}[ht]
\centering
\subfloat[without control]{\label{f5a}\includegraphics[width=0.31\textwidth]{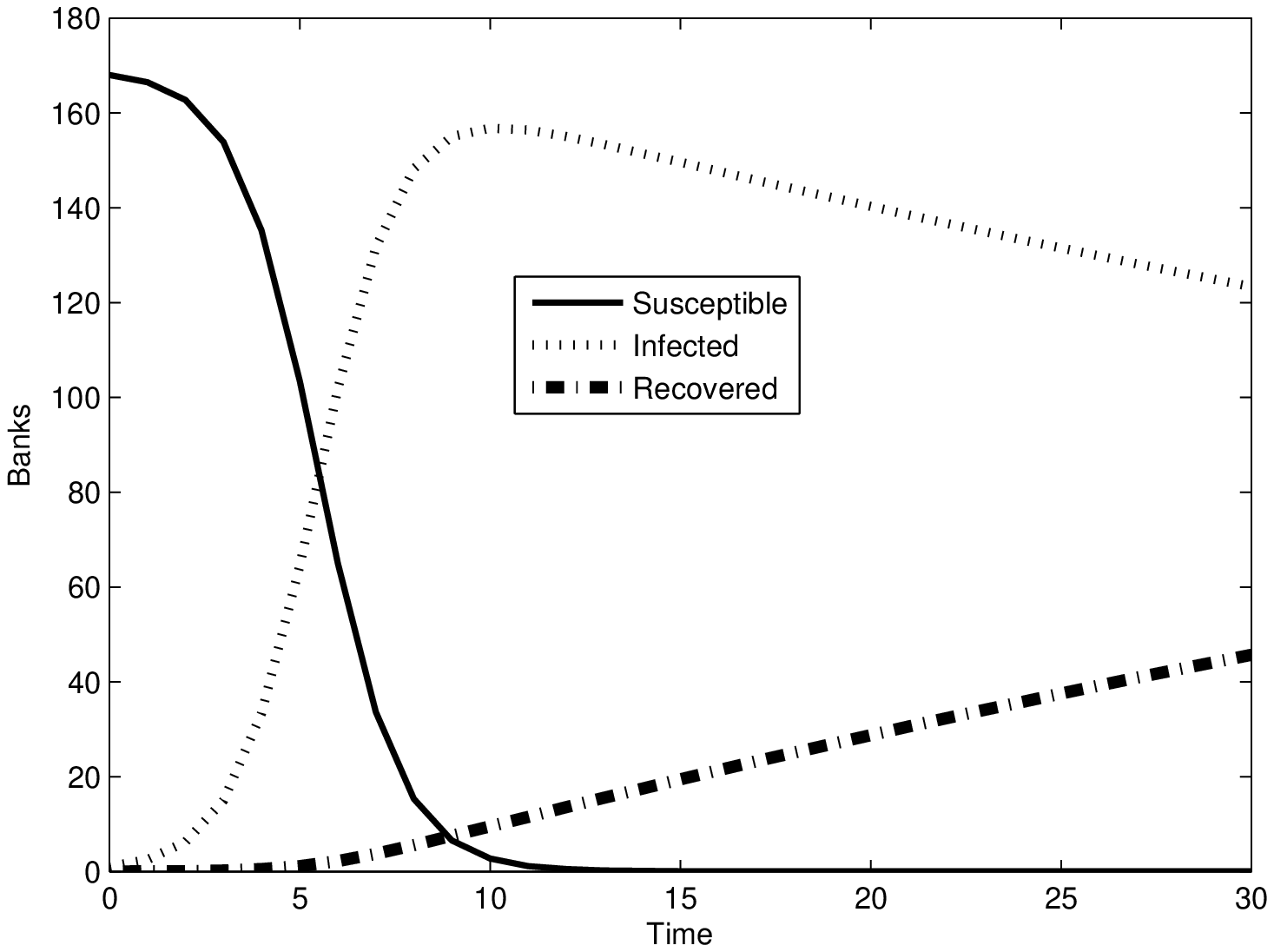}}
\subfloat[with optimal control]{\label{f5b}\includegraphics[width=0.33\textwidth]{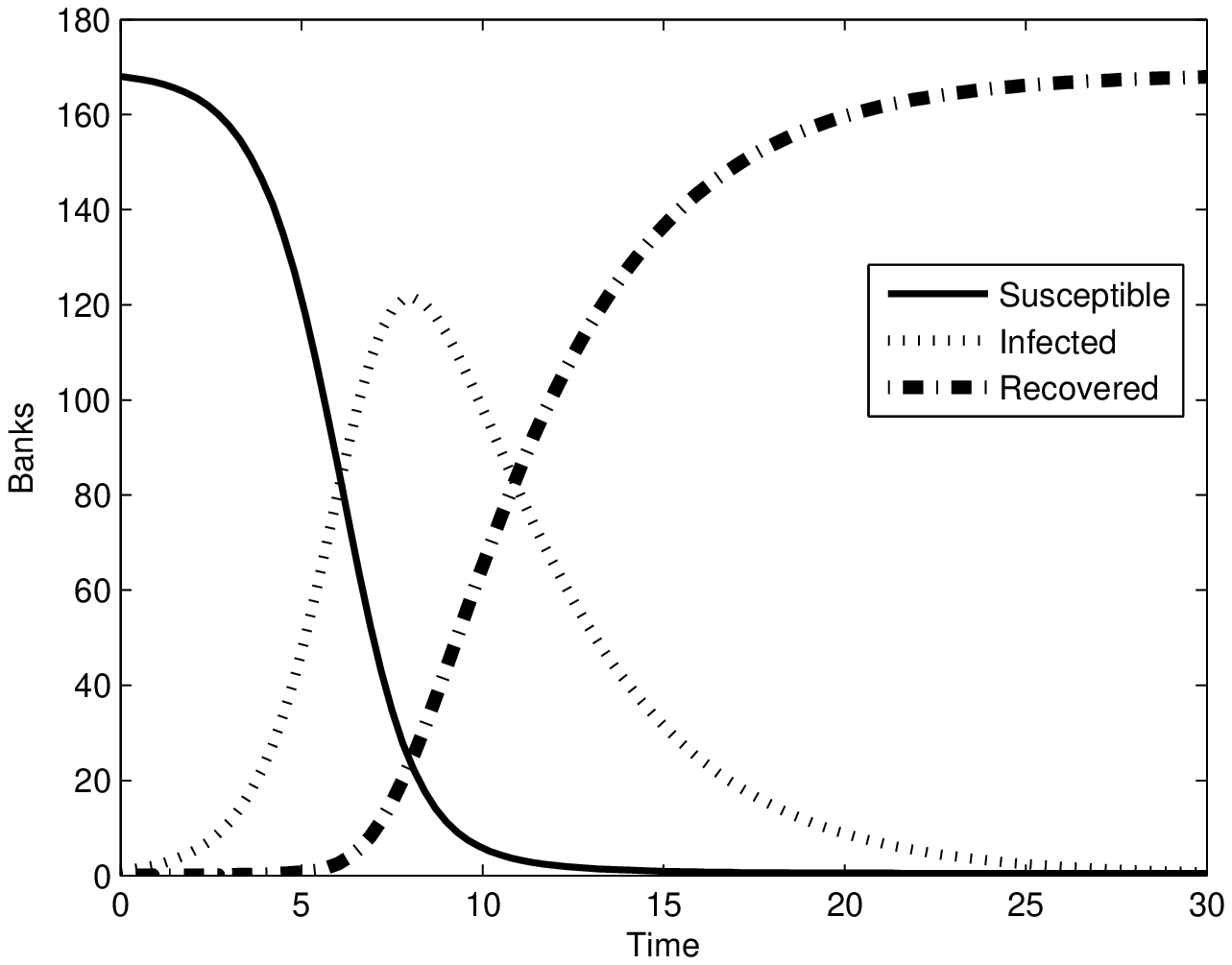}}
\subfloat[the extremal control $u(t)$]{\label{f5c}\includegraphics[width=0.33\textwidth]{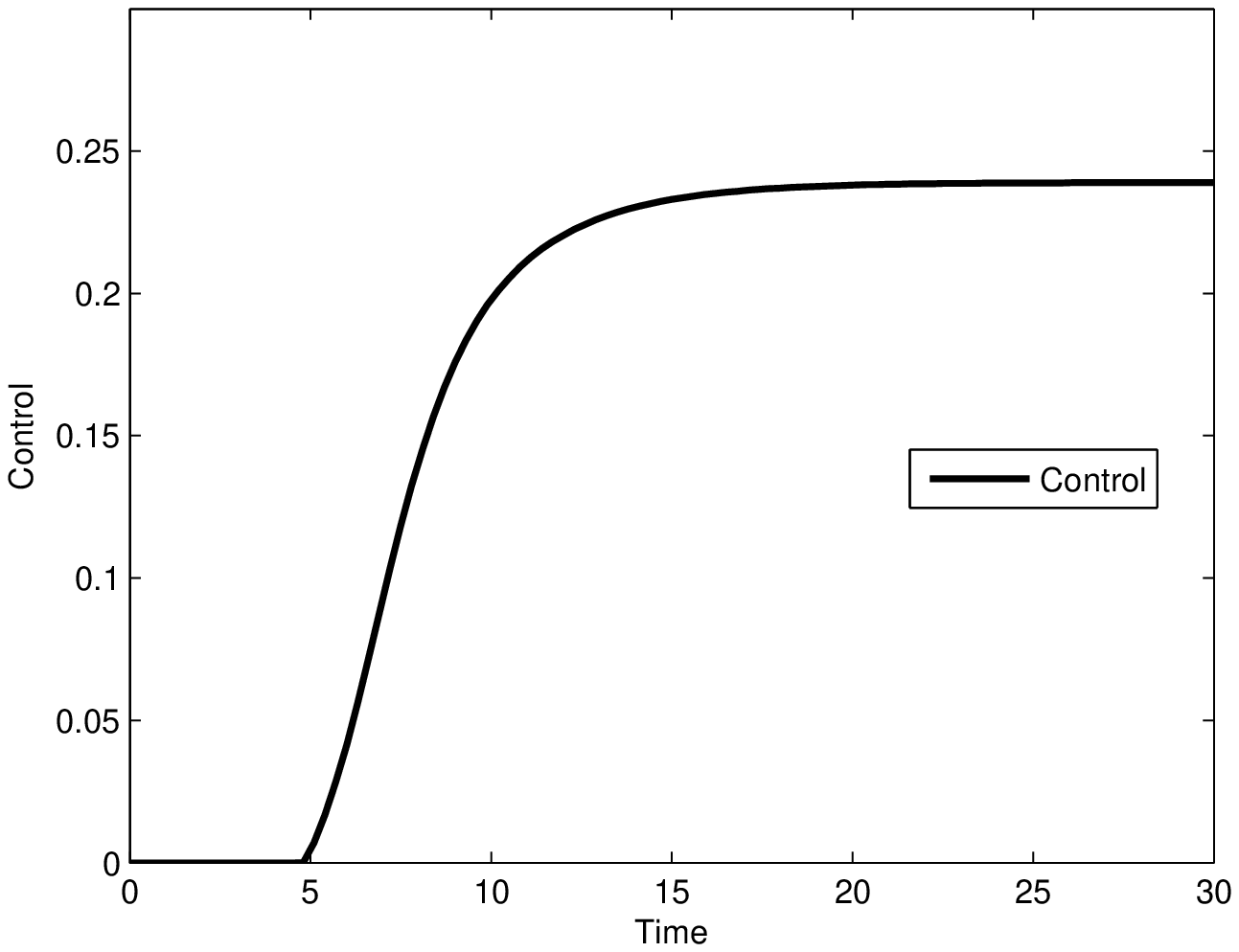}}
\vspace*{0.3cm}
\caption{Contagion risk from Spain (Scenario~2) with and without control.}
\label{f5}
\end{figure}

\begin{figure}[ht!]
\centering
\subfloat[without control]{\label{f6a}\includegraphics[width=0.33\textwidth]{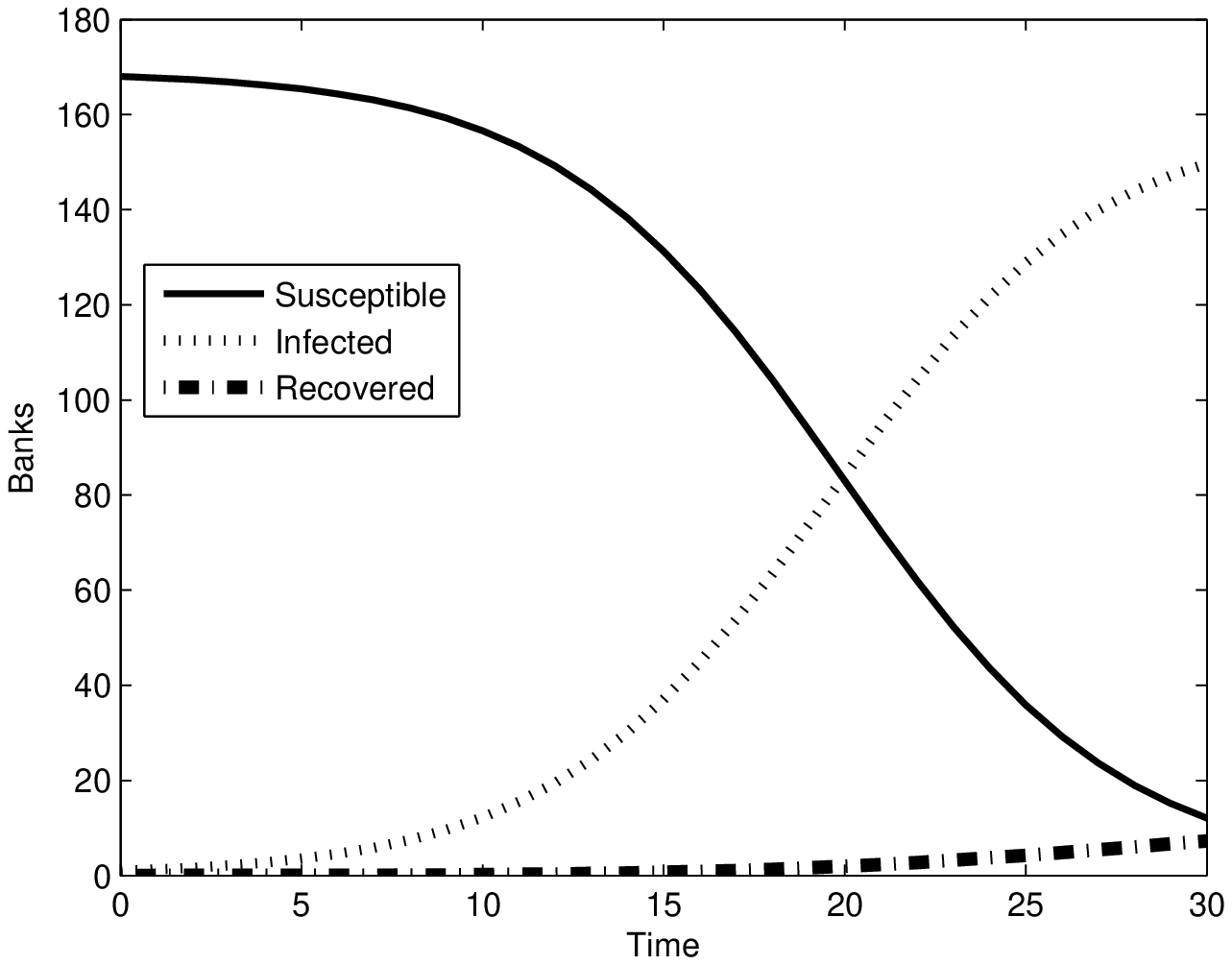}}
\subfloat[with optimal control]{\label{f6b}\includegraphics[width=0.33\textwidth]{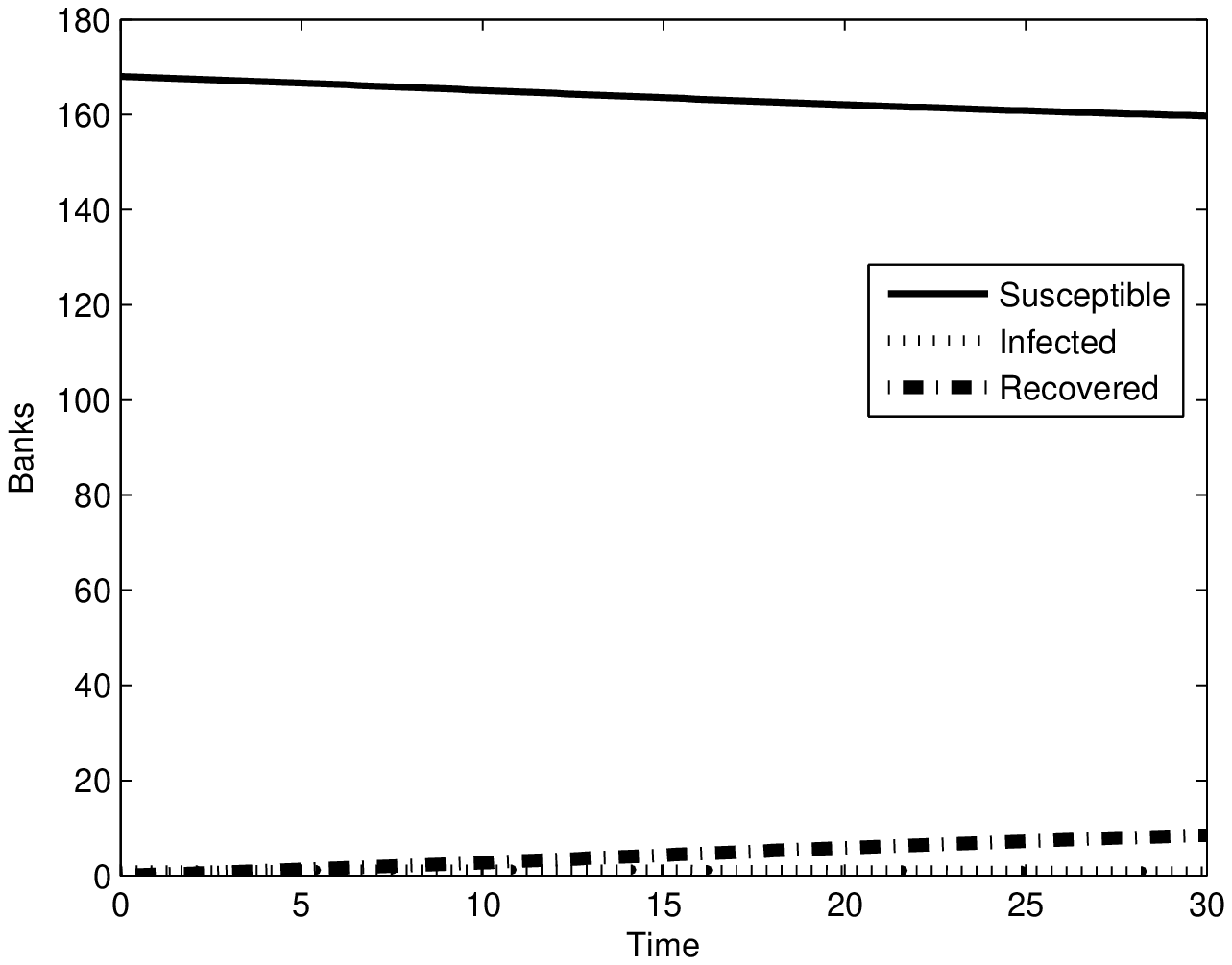}}
\subfloat[the extremal control $u(t)$]{\label{f6c}\includegraphics[width=0.33\textwidth]{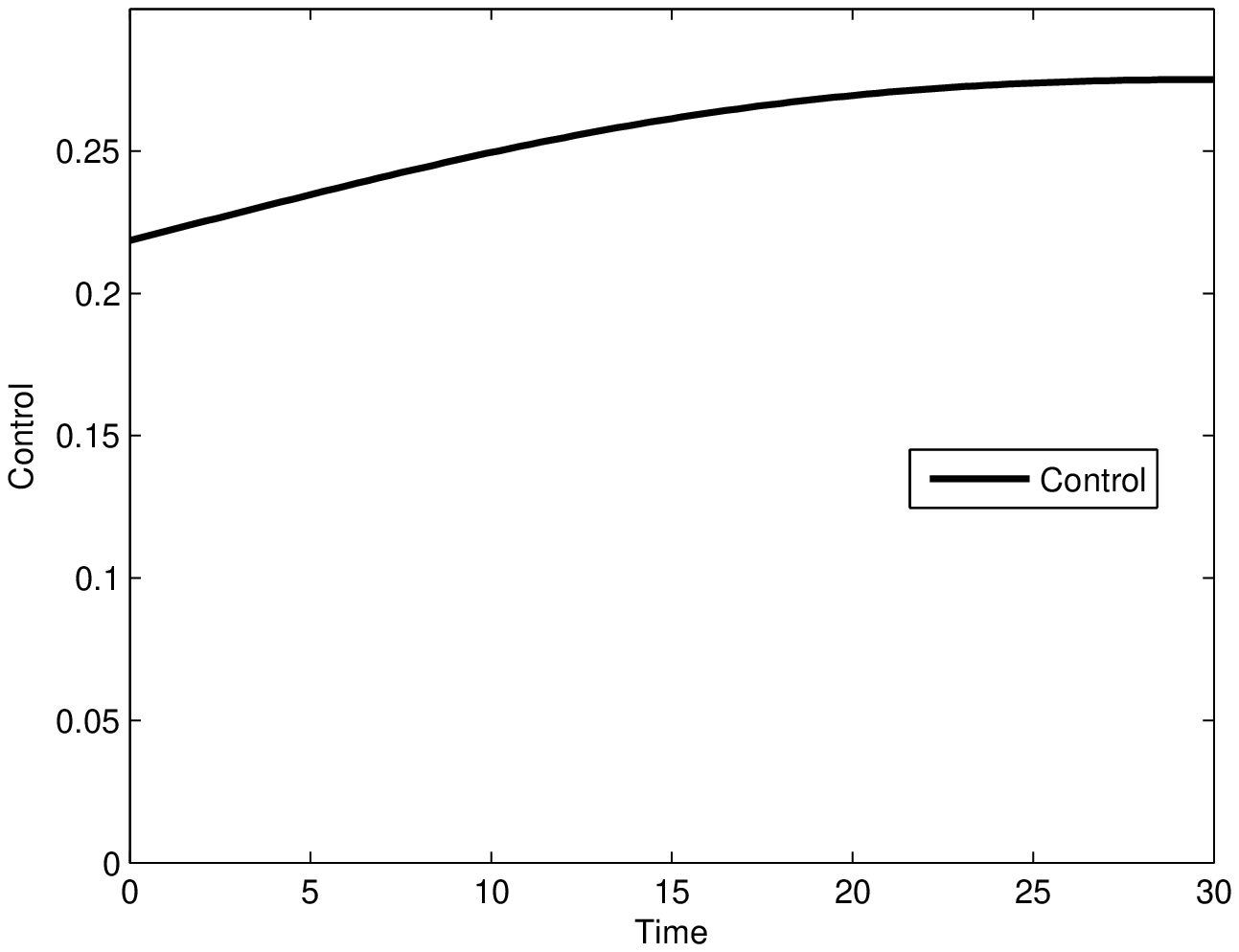}}
\vspace*{0.3cm}
\caption{Contagion risk from UK (Scenario~3) with and without control.}
\label{f6}
\end{figure}

In order to use the \textsf{BOCOP} optimal control solver \cite{[25],Bocop},
the optimal control problem in Bolza form \eqref{eqopt1}--\eqref{eqopt2} 
is rewritten in the following equivalent Mayer form:
\begin{equation}
\label{eqbocop1}
\min \mathcal{J}[S(\cdot),R(\cdot),Y(\cdot)]
=N-S(T)-R(T)+Y(T)\longrightarrow \min
\end{equation}
subject to
\begin{equation}
\label{eqbocop2}
\begin{cases}
\displaystyle \frac{dS(t)}{dt} = \beta S^2(t)+ \beta S(t)(R(t)-N),\\[0.3cm]
\displaystyle \frac{dR(t)}{dt} = \gamma (N-S(t)-R(t))+u(t)(N-S(t)-R(t)),\\[0.3cm]
\displaystyle \frac{dY(t)}{dt} = bu^2(t),\\[0.3cm]
\displaystyle S(0)=S_0, \quad R(0)=0, \quad Y(0)=0,\\[0.3cm]
u(t) \in [0, 1].
\end{cases}
\end{equation}
The optimal control problem \eqref{eqbocop1}--\eqref{eqbocop2} is approximated 
by \textsf{BOCOP} into a nonlinear programming problem (NLP) using a standard 
time discretization. The NLP problem is then solved by the open-source
\textsf{Ipopt} optimization solver \cite{ipopt}, using sparse exact derivatives 
computed by the \textsf{Adol-C} (Automatic Differentiation by OverLoading in \textsf{C++}) 
package \cite{ADOL-C}. Figures~\ref{f4}, \ref{f5} and \ref{f6} show the results 
of simulations of our SIR bank contagion risk model, with and without optimal control, 
for a period of 30 days ($T=30$) for the first, second and third scenarios 
of Section~\ref{s3}, respectively.

Figure~\ref{f4a} shows that if we have no interest to stop contagiousness 
or just do not want to spend any money on it, that is, $u(t) \equiv 0$, 
then, in case of the first scenario, the number of contagious banks at final 
time $T=30$ will be equal to 64. If the control $u(t)$ is chosen as the
solution of our optimal control problem, then in one month the number 
of contagious banks will become smaller than $\mathcal{J}=2$, see Figure~\ref{f4b}.

In case of the second scenario, Figure~\ref{f5a} shows that the number of contagious 
without control is equal to 124. In contrast, the number of contagious banks at time $T=30$
using optimal control, taking into account the costs associated with the control interventions
during all the period $[0,30]$, is equal to 2 (see Figure~\ref{f5b}). 
On the other hand, after one month, the number of banks 
that already recovered from the risks are less than fifty without control; 
while with optimal control measures this number is almost double. 

Finally, when the contagion begins from United Kingdom (third scenario), 
the number of contagious banks at final time $T=30$ is 149 without control 
(Figure~\ref{f6a}) with that number decreasing to less than 4 when the Central Bank 
applies optimal control (see Figure~\ref{f6b}).

The results of the three scenarios just reported, with respect to the value
of the cost functional $\mathcal{J}$, are summarized in Table~\ref{table1}
for a better comparison.
\begin{table}[ht!]
\caption{Values of the cost functional $\mathcal{J}$ \eqref{eqopt1}.}
\begin{tabulary}{\linewidth}{lcc} \hline
Scenario \# & no control measures & with optimal control\\ \hline
Scenario 1 \footnotesize{(bank contagion starting in Portugal)} & 64 & 2 \\
Scenario 2 \footnotesize{(bank contagion starting in Spain)} & 124 & 2  \\
Scenario 3 \footnotesize{(bank contagion starting in UK)} & 149 & 4 \\ \hline
\end{tabulary}
\label{table1}
\end{table}

\begin{figure}[ht]
\centering
\subfloat[Portugal]{\label{f7a}\includegraphics[width=0.33\textwidth]{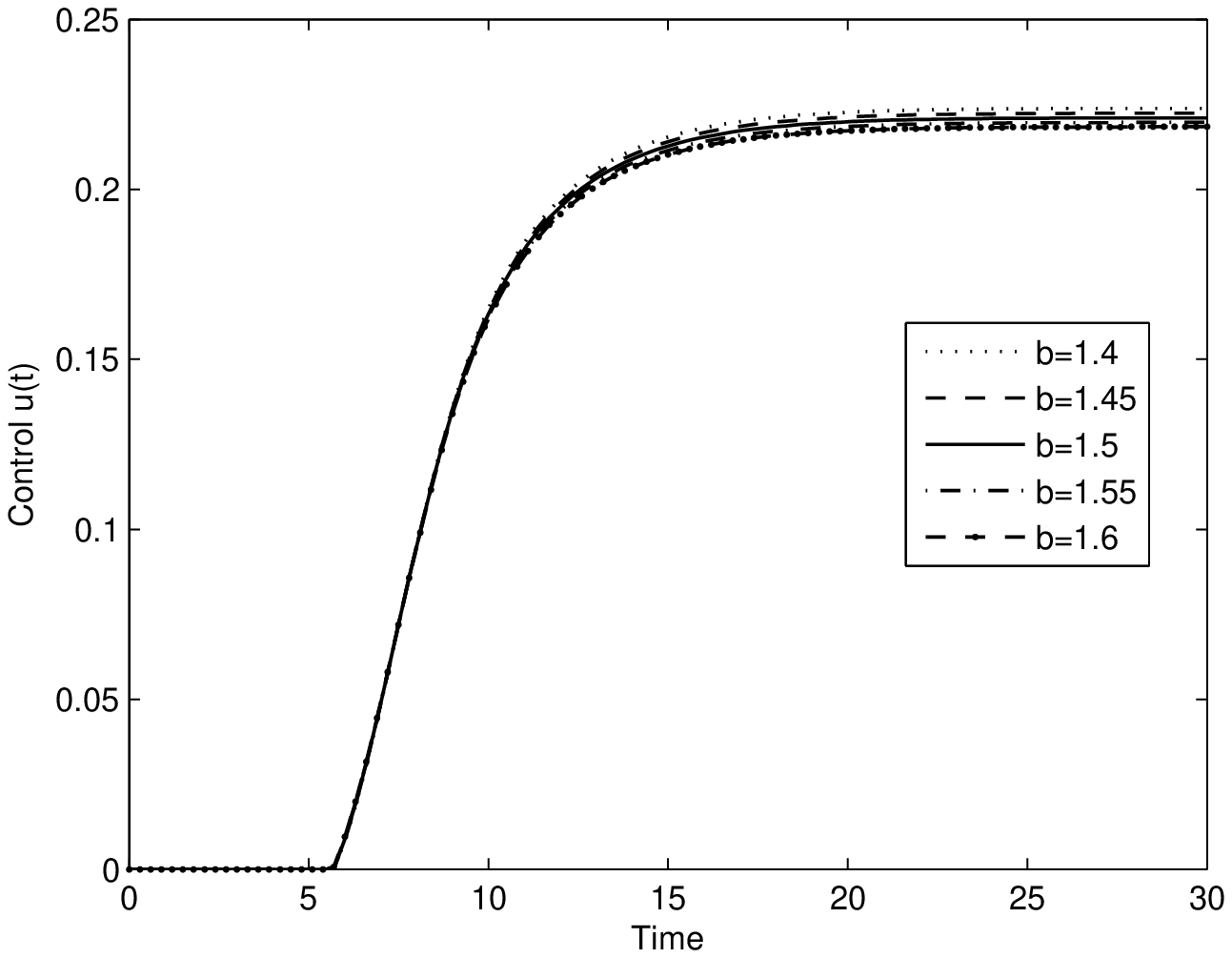}}
\subfloat[Spain]{\label{f7b}\includegraphics[width=0.33\textwidth]{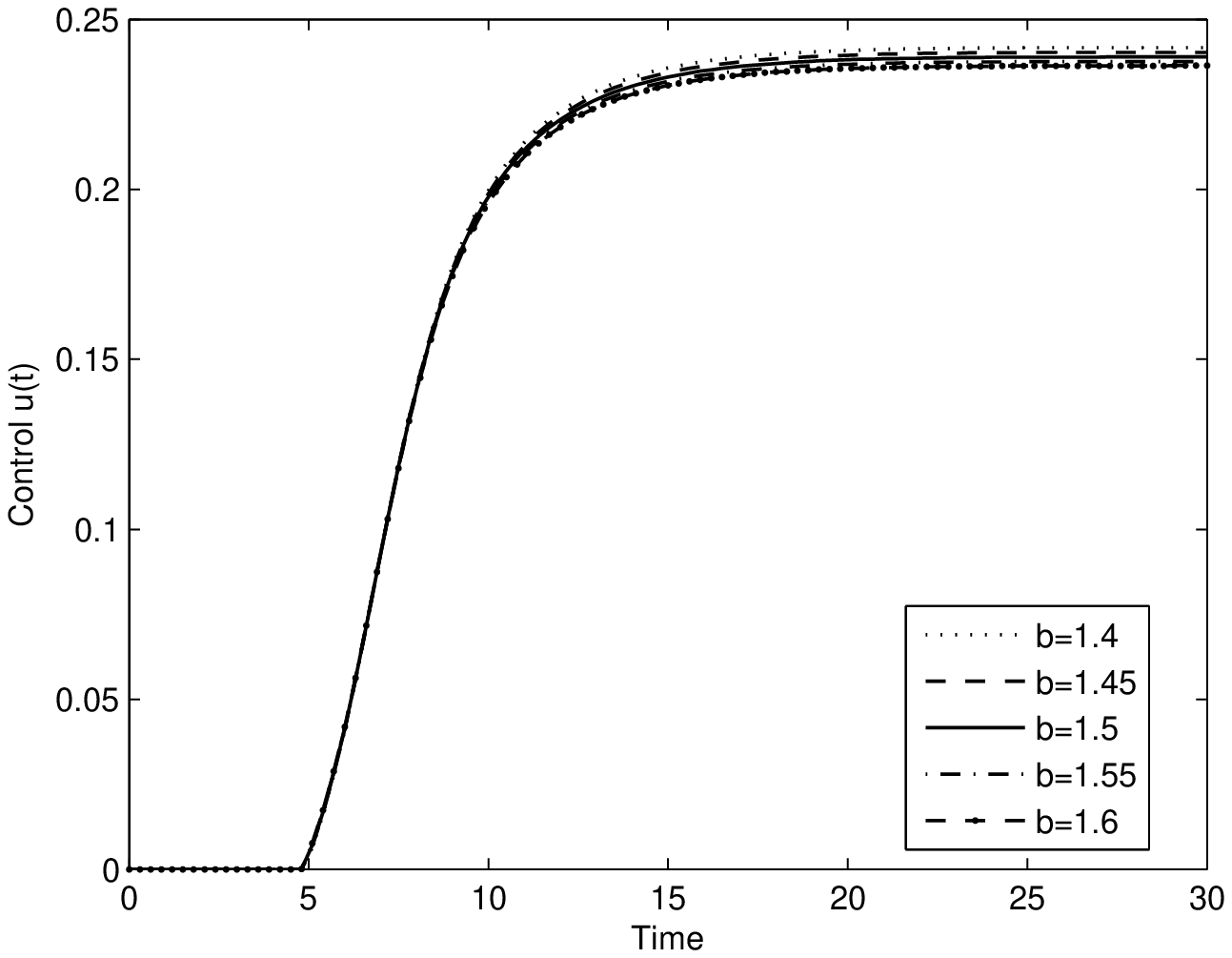}}
\subfloat[UK]{\label{f7c}\includegraphics[width=0.33\textwidth]{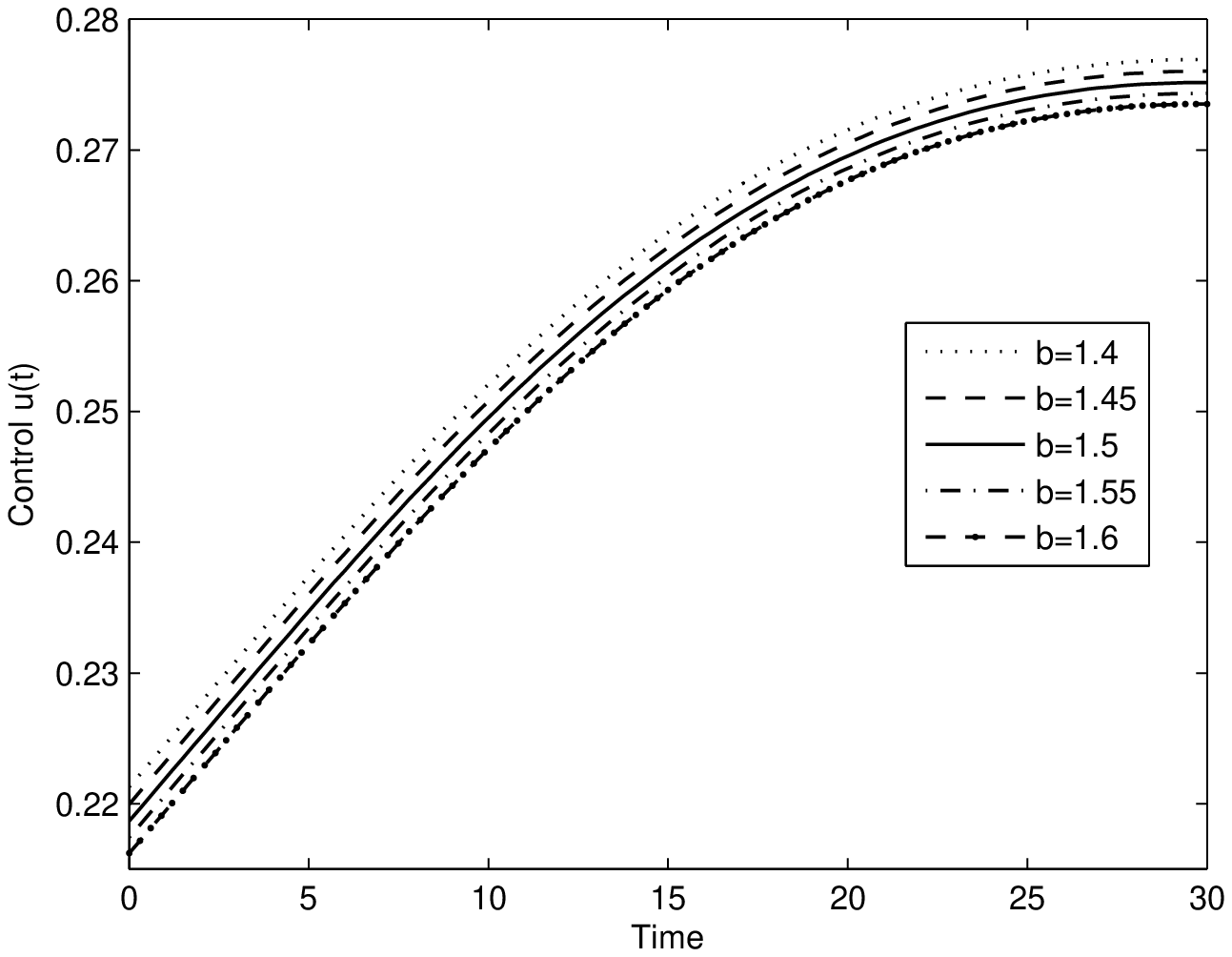}}
\vspace*{0.3cm}
\caption{The optimal control $u(t)$ for different values of the weight $b$.}
\label{f7}
\end{figure}

It is worth to mention that the qualitative results of the paper
do not change with the particular choice of the weight $b$: 
the extremal curves change continuously with the change of parameter $b$. 
For illustrative purposes, we show in Figure~\ref{f7}
the optimal control $u(t)$, $t \in [0, 30]$ days, 
for different values of the weight $b$.
We see that an increase of the value of the parameter $b$ 
corresponds to a decrease of the control $u$,
meaning that the financial support from the Central Bank 
diminishes when we penalize the use of the control by increasing the value of $b$.
The role of parameter $b$ is more visible in Scenario~3 than in Scenarios 1 and 2.

\section{Conclusions}
\label{s5}

We investigated the dynamic behaviour of contagiousness in the banking sector, 
at the macroeconomic level, using a SIR epidemic model and optimal control.
The features of contamination were shown to depend on the parameter values 
of the transmission and recovery rates, as well as on the country for which 
the process of infection begins. The scale of negative consequences 
for three different scenarios of bank risk contagion were identified. 
Banks at risk from countries with greater financial influence in Europe 
tend to propagate more severely the banking crisis; at the same time, 
the recovery period is longer. An optimal control problem was proposed, 
in order to reduce the number of contagious banks, to prevent large-scale 
epidemic contagiousness, and to avoid serious financial and economic consequences.

		
\section*{Acknowledgements}

This research was supported in part by the Portuguese Foundation 
for Science and Technology (FCT -- Funda\c{c}\~{a}o para a Ci\^{e}ncia e a Tecnologia), 
through CIDMA -- Center for Research and Development in Mathematics and Applications, 
within project UID/MAT/04106/2013. Kostylenko is also supported 
by the Ph.D. fellowship PD/BD/114188/2016.

We are very grateful to the authors of \cite{[24]}
for providing us the parameter values that they have obtained
in their work; and to two anonymous referees 
for valuable remarks and comments, which significantly 
contributed to the quality of the paper.
	


\end{document}